\documentclass[11pt]{article}

\usepackage{amssymb}
\usepackage{amsmath}
\usepackage{float}
\usepackage{graphicx}
\usepackage[none]{hyphenat}
\usepackage[english,american]{babel}

\newtheorem{theo}{Theorem}[section]
\newtheorem{defn}[theo]{Definition}
\newtheorem{nota}[theo]{Notation}
\newtheorem{exa}[theo]{Example}
\newtheorem{rema}[theo]{Remark}
\newtheorem{coro}[theo]{Corollary} 
\newtheorem{lemm}[theo]{Lemma}
\newtheorem{prop}[theo]{Proposition}

\newtheorem{ass}[theo]{Assumption}

\newcommand{\bigslant}[2]{{\raisebox{.2em}{$#1$}\left/\raisebox{-.2em}{$#2$}\right.}}
\def\R{\mathbb{R}}

\def\N{\mathbb{N}}
\def\Z{\mathbb{Z}}
\def\qed{\hfill\ensuremath\square}

\allowdisplaybreaks

\begin{document}

\title{The Jordan-H\"older Theorem for Monoids with Group Action}

\author{Alfilgen Sebandal\\
	{\small E-Mail:alfilgen.sebandal@g.msuiit.edu.ph}\\[.3cm]
	Jocelyn Vilela\\
	{\small E-Mail:jocelyn.vilela@g.msuiit.edu.ph}\\~\\
	{\small Center of Graph Theory, Algebra and Analysis}\\
		{\small	Premier Research Institute of Science and Mathematics}\\
			{\small Department of Mathematics and Statistics}\\
			{\small College of Science and Mathematics}\\
			{\small MSU-Iligan Institute of Technology
			Tibanga, Iligan City, Philippines}}

\maketitle

\begin{abstract}
	In this article, we prove an isomorphism theorem for the case of refinement $\Gamma$-monoids. Based on this we show a version of the well-known Jordan-H\"older theorem in this framework. The central result of this article states that - as in the case of modules - a monoid $T$ has a $\Gamma$-composition series if and only if it is both $\Gamma$-Noetherian and $\Gamma$-Artinian. As in module theory, these two concepts can be defined via ascending and descending chains respectively.
\end{abstract}
\vspace{0.3cm}


\section{Introduction}
The theorem of Jordan-H\"older was extended many times since its original proof in 1870 by C.~Jordan for groups \cite{Jor69, Jor70, Hol89, Sch28}. Based on this, one obtains isomorphism uniqueness for a composition series of simple factor groups. Nowadays the concept of composition series plays a key role in the analysis of several algebraic structures as groups, modules, algebras or categories see e.g.~the mono\-graphs \cite{Erdmann, Kashiwara} and the articles \cite{Natale, Cassidy, Natale2, Orlik} and the references therein. Especially for non-semi-simple modules, a composition series, that is a finite increasing chain of simple sub-constituents replaces the direct sum of simple ones and hence gives a profound access to the algebraic structure. While for many groups, modules and algebras the above theory is well-established and used in the corresponding representation theory, for monoids the contrary is the case. This is due to the fact that the corresponding isomorphism theorems have to be proven to hold, see e.g.~\cite{Lederman, Fritsch} for conditions under which a Jordan-H\"older theorem holds.  In this article we extend these concepts to abelian refinement $\Gamma$-monoids, this means monoids with a group acting on them. The concept of $\Gamma$-monoids and $\Gamma$-order-ideals of such monoids was recently systematically introduced in \cite{A1} to describe a monoid associated to graph algebras. The article is organized as follows: We first list the notations and definitions needed throughout the article. Then we give properties on $\Gamma$-monoids. Moreover, we prove the necessary isomorphism theorems and the Jordan-H\"older theorem. With the help of this we can show that a monoid T has a $\Gamma$-composition series if and only if it is both $\Gamma$-Noetherian and $\Gamma$-Artinian. As in module theory these two concepts can be defined via ascending and descending chains, respectively.

\section{Preliminaries}
The following are found in \cite{monoid.q} and \cite{A1}. A monoid is a set equipped with an associative binary operation and an identity element. 
Due to the lack of an inverse element, which allows for cancellations, the study of monoids is performed with the help of equivalence relations. Throughout this article we use two of these binary relations, which help us to classify ideals later on.

\begin{defn}
	Let $M$, $M_1$ and $M_2$ be commutative monoids. 
	\begin{enumerate}
		\item[\textnormal{i)}] For any submonoid $H$ of $M$,  we define a binary relation $\rho_H$ in $M$ by 
		$x~\rho_H ~y$ if~and~only~if $(x+H)\cap (y+H)\neq \varnothing.$
		\item[\textnormal{ii)}]  For a mapping $f:M_1\rightarrow M_2$, define a relation 
		$x\rho_f y$ if and only if $f(x)=f(y).$
	\end{enumerate}
\end{defn}
It can be shown that both $\rho_f$ and $\rho_H$ are equivalence relations for any mapping $f$ and submonoid $H$. For any submonoid $H$ of $M$, the set
$$\bigslant{M}{H} \cong \bigslant{M}{\rho_H} =\lbrace \rho_{H}(x):x\in M\rbrace \hfill$$
is an abelian monoid under the  operation $\circ$ defined by $\rho_H(x)\circ \rho_H(y) := \rho_H(x+y)$ with $\rho_H(0)$ as its identity and where $0$ is the identity in $M$. Furthermore, we say that $H$ is $normal$ if for any $x,y\in M$, $x,x+y\in H$ implies $y\in H$. Equivalently, $H$ is normal if $\rho_H(0)=H$.\\

There is a natural algebraic pre-ordering on the commutative monoid $M$ defined by $a\leq b$ if $b=a+c$, for some $c\in M$. Throughout, $a\parallel b$ shall mean the elements $a$ and $b$ are not comparable.
\begin{defn}
	A commutative monoid $M$ is called 
	\begin{itemize} 
		\item[\textnormal{i)}] $conical$, if $a+b=0$ implies $a=b=0$; 
		\item[\textnormal{ii)}] $cancellative$, if $a+b=a+c$ implies $b=c$, where $a,b,c\in M$; 
		\item[\textnormal{iii)}] $refinement$, if for $a+b=c+d$, there exist $e_1,e_2,e_3,e_4\in M$ such that $a=e_1+e_2$, $b=e_3+e_4$ and $c=e_1+e_3$, $d=e_2+e_4$. 
	\end{itemize}
	An element $a\in M$ is called 
	\begin{itemize} 
		\item[\textnormal{iv)}] $minimal$ if $b\leq a$ implies $a\leq b$. 
	\end{itemize}
\end{defn}

\begin{rema}If $M$ is conical and cancellative, these notions coincide with the more intuitive definition of minimality, this means $a$ is minimal if $0\neq b\leq a$, then $a=b$.
\end{rema}
Next we define the main objects needed throughout this work.

\begin{defn}
	\textnormal{An \emph{action} of a group $\Gamma$ on a set $M$ is a function $\Gamma\times M\rightarrow M$ ($(\alpha, a)\mapsto {^\alpha a}$)  such that for all $a\in M$ and $\alpha, \beta\in \Gamma$, $^0a=a$ and $^{\alpha \beta}a= {^\alpha}(^\beta a)$. Let $M$ be a monoid with a group $\Gamma$ acting on it. Then $M$ is said to be a  $\Gamma$-\emph{monoid}. For $a\in M$, denote the orbit of the action
		of $\Gamma$ on an element a by $O(a)$, so $O(a) = \{{^\alpha a} | \alpha \in \Gamma\}$.}
	
	An action of a group $\Gamma$ on a set $M$ with algebraic structure needless to say must be compatible to the operations on the set. Hence, for any $M$ monoid, we have $^\alpha (a+b)= {^\alpha}a+{^\alpha}b$.

	\begin{exa}
		Let $\Gamma$ be the set of integers $\Z$ under addition and let $T=M_2(\R)$ under pointwise addition. It could be shown that $\Z\times T\rightarrow T$ given by $$ \left ( x, \begin{pmatrix} a & b\\
		c & d\end{pmatrix}
		\right )\mapsto \begin{pmatrix} e^x a & b\\
		c & e^x d\end{pmatrix}$$ is an action which makes $T$ a $\Gamma$-monoid.
	\end{exa}

\end{defn}
The study of ideals is traditionally linked with homomorphisms. This is the same in our setting.
\begin{defn} \textnormal{Let $M, M_1$ and $M_2$ be monoids. $\Gamma$ a group acting on $M, M_1$ or $M_2$, respectively.
		\begin{itemize}
			\item[\textnormal{i)}] A $\Gamma$-\emph{module} \emph{homomorphism} is a monoid homomorphism $\phi:M_1\rightarrow M_2$ that respects the action of $\Gamma$, this means $\phi(^\alpha a)={^\alpha \phi(a)}$.
			\item[\textnormal{ii)}] A $\Gamma$-\emph{order-ideal} of a monoid $M$ is a subset $I$ of $M$ such that for any $\alpha, \beta\in \Gamma$, $^\alpha a+{^\beta b}\in I$ if and only if $a,b\in I$.
	\end{itemize}}
\end{defn}
\begin{rema}
	Equivalently, a $\Gamma$-order-ideal is a submonoid $I$ of $M$ which is closed under the action of $\Gamma$ and it is $hereditary$ in the sense that $a\leq b$ and $b\in I$ imply $a\in I$. 
\end{rema}

The set $\mathcal{L}(M)$ of $\Gamma$-order-ideals of $M$ forms a complete lattice. We say $M$ is a $simple$ $\Gamma$-$monoid$ if the only $\Gamma$-order-ideals of $M$ are $0$ and $M$.\\
\begin{nota} For $a\in M$, we denote the $\Gamma$-order-ideal generated by an element $a$ by $\langle a \rangle$. It is easy to see that 
	$$\langle a \rangle =\left \{  x\in M: x\leq \sum_{\alpha \in \Gamma}{^\alpha a} \right \}.$$
\end{nota}

\section{Some Properties of $\Gamma$-monoids}

\begin{ass}
	Throughout this paper, we shall assume that the group $\Gamma$ is commutative and we let $T$ be a $\Gamma$-monoid with identity denoted by $0$ and the operation $\circ$ on the quotient monoid by $+$.
\end{ass}

\begin{rema}\label{pro35}
	Every $\Gamma$-order-ideal of $T$ is normal.
\end{rema}

\begin{rema}\label{rem36} Let $I$ be a $\Gamma$-order-ideal of $T$. Then by the properties of an equivalence relation,
	\begin{itemize}
		\item[\textnormal{i)}] for every $g\in T $, $g\in \rho_I(g)$ and
		\item[\textnormal{ii)}] $\rho_I(g_1)=\rho_I(g_2)$ if and only if $(g_1+ I)\cap (g_2+I)\neq \varnothing$.
	\end{itemize}
\end{rema}

\begin{coro}\label{cor36b} Let $I$ be a $\Gamma$-order-ideal of $T$. Then  $\rho_I(g)=\rho_I(0)$ if and only if $g\in I$.
\end{coro}

\noindent \textbf{Proof:}{ Suppose $g\in I$. Let $x\in \rho_I(g)$. Then $(x+I)\cap (g+I)\neq \varnothing$. Thus, there exist $h_1,h_2\in I$ such that $x+h_1=g+h_2\in I$.
	Since $I$ is a $\Gamma$-order-ideal, we have $x\in I$. Thus, $x+0=0+x$ implies 
	$(x+I)\cap (0+I)\neq \varnothing. $
	Hence, $x\in \rho_I(0)$. It follows that $\rho_I(g)\subseteq \rho_I(0)$. \\
	\indent Let $x\in \rho_I(0)$. Then $(x+I)\cap (0+I)\neq \varnothing$. Thus, there exist $h_1,h_2\in I$ such that 
	$x+h_1=0+h_2\in I.$ 
	Since $I$ is a $\Gamma$-order-ideal, $x\in I$. Now, since $g\in I$, $x+g=g+x$ implies $(x+I)\cap (g+I)\neq \varnothing$. Thus, $x\in \rho_I(g)$. Thus, $\rho_I(0)\subseteq \rho_I(g)$.\\
	\indent Accordingly, $\rho_I(0)= \rho_I(g)$.\\
	\indent Now, suppose, $\rho_I(0)= \rho_I(g)$. Then by Remark \ref{rem36}(ii), we have $(g+I)\cap (0+I)\neq \varnothing$. Thus, there exist $h_1,h_2\in I$ such that  
	$g+h_1=0+h_2=h_2\in H.$ Since $I$ is a $\Gamma$-order-ideal, it follows that $g\in I$.~~ \qed
}

\begin{prop}\label{Prop36}
	An action of $\Gamma$ on $T$  induces  an action of $\Gamma$ on $T/I$ for any $\Gamma$-order ideal $I$ of $T$.
\end{prop}
\noindent \textbf{Proof:} Define a mapping $\Gamma\times T/I\rightarrow T/I$ by $^\alpha \rho_I(x)= \rho_I(^\alpha x)$ for all $\alpha \in \Gamma$ and $x\in T$. Let $(\alpha, \rho_I(x)), (\beta, \rho_I(y))\in \Gamma \times T$ such that $(\alpha, \rho_I(x))= (\beta, \rho_I(y))$. Then $\alpha =\beta$ and $\rho_I(x)=\rho_I(y)$. Thus, $(x+I)\cap (y+I)\neq \varnothing$ which implies that $x+z_1=y+z_2$, $z_1,z_2\in I$. 
Accordingly, 
$$^\alpha x + {^\alpha z_1}={^\alpha} (x+z_1)={^\beta } (y +z_2)={^\beta y} +{^\beta z_2}.$$
Since $^\alpha z_1, {^\beta} z_2 \in I$, being a $\Gamma$-order ideal,  $({^\alpha}x+I)\cap ({^\beta}y+I)\neq \varnothing$, this means 
$$^\alpha \rho_I(x)= \rho_I(^\alpha x)=\rho_I({^\beta y})={^\beta}\rho_I(y).$$
This implies well-definedness. \qed

\begin{prop}Let $T$ be a monoid and $I$ a submonoid of $T$ such that a group $\Gamma$ acts on $T/I$. If for all distinct elements $x,y\in T$, we have $(x+I)\cap (y+I)=\varnothing$,  then the action on $T/I$  by $\Gamma$ induces an action on $T$.
\end{prop}

\noindent \textbf{Proof:} Define a mapping $\Gamma\times T\rightarrow T$ by $(\alpha,x)\mapsto {^\alpha x}=y$ where $y\in T$ such that $^\alpha \rho_I(x)=\rho_I(y)$.\\ \indent Let $(\alpha,x)=(\beta, z)$, $\alpha, \beta \in \Gamma$ and $x,z\in T$.  Then $\alpha=\beta$ and $x=z$. Thus, $^\alpha \rho_I(x)={^\beta \rho_I(z)}\in T/I.$
This implies that 
$^\alpha \rho_I(x)=\rho_I(y)={^\beta \rho_I(z)}$ for some $y\in T$. By assumption, $(y+I)\cap (w+I)=\varnothing$, that is,  $\rho_I(y)\neq \rho_I(w)$ for every $w\neq y$. Hence, $y$ is unique of such. Thus, in $T$, $^\alpha x=y ={^\beta z}$. This implies well-definedness of the mapping on $T$.\\
\indent Now, for every $x\in T$,  $^0\rho_I(x)=\rho_I(x)$. Thus, $^0x=x$ in $T$. \\
\indent Let $\alpha, \beta \in \Gamma$ and $x\in T$. We show that $^{(\alpha +\beta)} x={^\alpha} ({^\beta}x)$. Now, $^{(\alpha +\beta)} x=z$ for some $z\in T$. Thus, $^{(\alpha +\beta)}\rho_I( x)= \rho_I(z)$. By the action in $T/I$, we have $^\alpha({^\beta \rho_I(x)})={^{(\alpha +\beta)}\rho_I( x)}=\rho_I(z   )$. Now, $^\beta x=t$ for some $t\in T$. Thus, $^\beta \rho_I(x)=\rho_I(t)$. Accordingly, $\rho_I(z   )={^{\alpha}( ^\beta\rho_I( x)})={ ^\alpha \rho_I(t) }$ which implies that $^\alpha t=z$. Hence, $z={^\alpha t}={^\alpha ({^\beta x})}$. Consequently,  $^{(\alpha +\beta)} x=z={^\alpha} ({^\beta}x)$.\\
\indent Therefore, $\Gamma$ acts on $T$.~ \qed

\begin{prop} \label{pro38}Let $I$ be a $\Gamma$-order-ideal of $T$. Then $T=I$ if and only if $T/I=\{ \rho_I(0) \}$.
\end{prop}

\noindent \textbf{Proof:} Suppose $T=I$. Let $x\in T/I=T/T$. Then $x=\rho_T(y)$ for some $y\in T$. By Corollary \ref{cor36b}, we have $\rho_T(0)=\rho_T(y)=x$. Hence, $T/T=T/I=\{\rho_T(0)  \}$. Conversely, suppose $T/I=\{\rho_I(0)  \}$. Let $x\in T$. Then $\rho_I(x)\in T/I$. Thus, $\rho_I(x)=\rho_I(0)$. By Corollary \ref{cor36b}, it follows that $x\in I$. Hence, $T\subseteq I$. Accordingly, $T=I$.~~\qed

\begin{lemm}\label{lemc3}Let  $A$ and $B$ be $\Gamma$-order-ideals of a refinement $\Gamma$-monoid $T$. Then $A+B$ is a $\Gamma$-order-ideal of $T$. 
\end{lemm}
\noindent \textbf{Proof:} Let $x,y\in A+B$ and $\alpha, \beta \in \Gamma$. Then $x=a_1+b_1$ and $y=a_2+b_2$ for some $a_1,a_2\in A$ and $b_1, b_2\in B$. Since $A$ and $B$ are $\Gamma$-order-ideals, $^\alpha a_1 +{^\beta a_2}\in A$ and  $^\alpha b_1 +{^\beta b_2}\in B$. Now, 
$$ {^\alpha x} +{^\beta y}= 
{^\alpha (a_1+b_1)} +{^\beta (a_2+b_2)}=(^\alpha a_1 +{^\beta a_2}) + (^\alpha b_1 +{^\beta b_2})\in A+B.$$ 

Conversely, suppose ${^\alpha x} +{^\beta y}\in A+B$ for every $\alpha, \beta\in \Gamma$. Then for $\alpha =0=\beta$, we have $^0x+{^0y}=x+y\in A+B$. Thus, $x+y=a+b$ for some $a\in A $ and $b\in B$. Since $T$ is a refinement monoid, we have $x=e_1+e_2$, $y=e_3+e_4$, $a=e_1+e_3$ and $b=e_2+e_4$ for some $e_1,e_2,e_3,e_4\in T$. Thus, $e_1+e_3=a\in A$ and $e_2+e_4=b\in B$. Since $A$ and $B$ are $\Gamma$-order-ideals, 
$e_1,e_3\in A$ and $e_2,e_4\in B$. Accordingly, 
$x=e_1+e_2\in A+B$ and  $y=e_3+e_4\in A+B$. \\
\indent Therefore, $A+B$ is a $\Gamma$-order-ideal of $T$. \qed

\begin{rema} Indeed the assumption of refinement monoid in Lemma \ref{lemc3} is crucial, since the set of order ideals form a lattice. If $T$ is not a refinement monoid, in general the sum of $\Gamma$-order-ideals is not a $\Gamma$-order-ideal.
\end{rema}

\begin{exa}Consider the set $T=\{0,1,x,y,z,s,b\}$ and an operation $(+)$ on given by
	$$
	\begin{array}{l| l l l l l l l}
	+ & 0& 1& x& y& z& s& b\\
	\hline
	0& 0 &1 &x &y &z &s &b\\
	1& 1 & 1 & 1 & s & s & s & b\\
	x & x & 1 & 1 & s & s & s & b\\
	y & y & s & s & y & y & s & b\\
	z & z & s & s & y & y & s & b\\
	s & s & s & s & s & s & s & b \\
	b & b & b & b & b & b & b & s
	\end{array} 
	$$
	Clearly the operation is commutative. For associativity and the closedness of subsets with respect to $+$, we carry out a more detailed computation.
	Let $A =\{0,1,x\}$ and $B =\{0,y,z\}$. It is easy to verify that $+$ is associative on $A$ and $B$. Now for all $u,u' \in A$ and $v,v' \in B$, such that $u,v \neq 0$, we have
	$$u+(v+v')=u+y=s=s+v=(u+v) +v'$$
	and	  
	$$
	(u+u')+v=1+v=s=u+s=u+(u'+v).
	$$
	Accordingly, $+$ is associative on $A\cup B$. Now, since for all $u,v \in A\cup B \cup \{s\}$
	$$
	(u+v)+s=s=u+s=u+(v+s),
	$$
	we have that $A\cup B \cup \{s\}$ is associative with respect to $+$. Furthermore for $u,u' \in A$, $v,v'\in B$
	\begin{align*}
		(u+v)+b &= s+b = b = u+b = u+(v+b)\\
		(u+u')+b &= s+b = b = u+b = u+(u'+b)\\
		(v+v')+b &= s+b = b = v+b = u+(v'+b).
	\end{align*}
	Hence we are left to show associativity for sums with $b$ and $s$ as summands. For $u \in A\cup B$,
	\begin{align*}
		(u+s)+b &= s+b = b = u+b = u+(s+b)\\
		(s+s)+b &= s+b = b = s+b = s+(s+b)\\
		s+(b+b) &= s+s = s = b+b = (s+b)+b. 
	\end{align*}
	Due to commutativity the associativity of $+$ is shown. Hence indeed $T$ is a monoid and $A$ and $B$ are submonoids of $T$. \\
	With $\Gamma=\{0\}$ acting trivially on $T$, we obtain that $T$ is a $\Gamma$-monoid. It can be verified easily that $A$ and $B$ are $\Gamma$-order ideals of $T$. \\
	Since $1+1=x+x$ can not be refined, $T$ is not a refinement monoid. 
	Moreover we have that
	$
	A+B = \{u+v:\, u\in A, \, v\in B\}= A\cup B\cup \{s\}.
	$   \\
	But we have
	$
	b+b =s = x+y \in A+B,
	$
	but $b\notin A+B$. Hence $A+B$ can not be a $\Gamma$-order ideal of $T$.
\end{exa}

\begin{lemm}\label{lemg1}Let $A$ and $B$ be $\Gamma$-order ideals of a refinement monoid $T$ such that $A\cap B=\{0\}$. Then 
	$\bigslant{(A+B)}{A}\cong B.$
\end{lemm}
\noindent \textbf{Proof:} Define a mapping $f:\bigslant{(A+B)}{A}\rightarrow  B$ by $f(\rho_A(a+b))=b$ for all $\rho_A(a+b)\in \bigslant{(A+B)}{A}$. \\
\indent Let $\rho_A(a+b),\rho_A(c+d)\in (A+B)/A$  such that $\rho_A(a+b)=\rho_A(c+d)$, where $a,c\in A$ and $b,d\in B$. By Corollary \ref{cor36b}, $\rho_A(a)=\rho_A(0)=\rho_A (c)$. Thus,
$\rho_A(b)=\rho_A(0)+\rho_A(b)=\rho_A(a)+\rho_A(b)=\rho_A(a+b).$
Similarly, $\rho_A(d)=\rho_A(c+d)$. Hence, 
$\rho_A(b)=\rho_A(a+b)=\rho_A(c+d)=\rho_A(d).$
By Remark \ref{rem36}(ii), there exist $x,y\in A$ such that $x+b=y+d$. Since $T$ is a refinement monoid, there exist $e_1,e_2,e_3,e_4\in T$ such that $x=e_1+e_2$, $b=e_3+e_4$, $y=e_1+e_3$ and $d=e_2+e_4$. 
Now, since $A$ and $B$ are $\Gamma$-monoids, $e_1+e_2=x,e_1+e_3=y\in A$  imply $e_1, e_2, e_3\in A$ and $e_3+e_4=b,e_2+e_4=d\in B$ implies $e_2,e_3, e_4\in B$. This implies that $e_2,e_3\in A\cap B=\{0 \}$. Thus, $x=e_1+e_2=e_1+0=e_1+e_3=y$
and $b=e_3+e_4=0+e_4=e_2+e_4=d.$
Thus, $f(\rho_A(a+b))=b=d=f(\rho_A(c+d))$. Hence, $f$ is well-defined. \\
\indent Now,
\begin{eqnarray*}
	f(\rho_A(a+b)+\rho_A(c+d)&=& f(\rho_A((a+c)+(b+d))\\
	&=& b+d\\
	&=& f(\rho_A(a+b))+f(\rho_A(c+d)).
\end{eqnarray*}
Thus, $f$ is a monoid homomorphism.\\
\indent Now suppose $f(\rho_A(a+b))=f(\rho_A(c+d))$. Then $b=d$. Thus, $\rho_A(b)=\rho_A(d)$. Since $a,c\in A$, $\rho_A(a)=\rho_A(0)=\rho_A(c)$. Accordingly, 
$$\rho_A(a+b)=\rho_A(a)+\rho_A(b)=\rho_A(c)+\rho_A(d)=\rho_A(c+d).$$
Hence, $f$ is one-one. \\
\indent Let $b\in B$. Then $0+b\in A+B$ and $f(\rho_A(0+b))=b$. Hence, $f$ is onto. \\
\indent Therefore, 
$\bigslant{(A+B)}{A}\cong B.$ \qed
~\\

\begin{lemm}\label{lem46}Let  $A$ and $B$ be $\Gamma$-order-ideals of $T$. Then also $A\, \cap\, B$ is a $\Gamma$-order-ideal of $T$.
\end{lemm}
\noindent \textbf{Proof:} Let $\alpha, \beta \in \Gamma$ and $x,y\in A\cap B $. Then $x,y\in A$ and $x,y\in B$. Since $A$ and $B$ are $\Gamma$-order-ideals of $T$, $^\alpha x+^\beta y \in A$ and $^\alpha x+^\beta y \in B$. Hence, $^\alpha x+^\beta y \in A\cap B$.\\
\indent Suppose $^\alpha x+^\beta y \in A\cap B$. Take $\alpha =0=\beta$. Then $x+y\in A\cap B$. Thus, $x+y\in A$ and $x+y\in B$. Since $A$ and $B$ are $\Gamma$-order-ideals of $T$, we have $x,y\in A$ and $x,y\in B$. Accordingly, $x,y\in A\cap B$.\\
\indent Therefore, $A\cap B$ is a $\Gamma$-order-ideal of $T$.~~\qed
~\\

Since a $\Gamma$-order ideal is a submonoid, $T/I$ is defined for any $\Gamma$-order ideal of $T$ and we have the following lemma. 
\begin{lemm}\label{lem57}Let  $I$ be a $\Gamma$-order-ideal of $T$. If $T/I$ is a $\Gamma$-monoid, then every $\Gamma$-order ideal of $T/I$ is of the form $J/I$ where $J$ is a $\Gamma$-order-ideal of $T$ containing $I$.\end{lemm}

\noindent \textbf{Proof:} Let $H$ be a $\Gamma$-order-ideal of $T/I$. Then $H\subseteq T/I$. Let $J=\{ g\in T:\rho_I(g)\in H \}$. We show that $J$ is a $\Gamma$-order-ideal of $T$.\\
\indent Let $x,y\in J$ and $\alpha, \beta \in \Gamma$. Then $\rho_I(x), \rho_I(y)\in H$ and  $^\alpha \rho_I(x) + {^\beta \rho_I(y)}\in H$, since $H$ is a $\Gamma$-order-ideal. Accordingly, we have
$$\rho_I(^\alpha x+ {^\beta y})= \rho_I(^\alpha x) +\rho_I(^\beta y)= {^\alpha \rho_I(x)} + {^\beta \rho_I(y)}\in H.$$
It follows that $^\alpha x+{^\beta y}\in J$. \\
\indent Conversely, suppose $^\alpha x+{^\beta y}\in J$ for every $\alpha , \beta\in \Gamma$. For $\alpha=0=\beta $, $x+y\in J$. It follows that $\rho_I(x) + \rho_I(y)=\rho_I(x+y)\in H$. Since $H$ is a $\Gamma$-order-ideal, we have $\rho_I(x), \rho_I(y)\in H$, that is, $x,y\in J$. Accordingly, $J$ is a $\Gamma$-order-ideal of $T$. Now, we show that $I\subseteq J$. Let $x\in I$. Then by Corollary \ref{cor36b}, we have $\rho_I(x)= \rho_I(0)$. Since $\rho_I(0)$ is the identity in $T/I$ and $H$ is a $\Gamma$-order-ideal of $T/I$, we must have $\rho_I(x)= \rho_I(0)\in H$. Thus, $x\in J$. Hence $I\subseteq J$.  \qed

\begin{theo}\textnormal{\cite{monoid.q}} \label{thm13}
	Let $F_1$ and $F_2$ be commutative monoids and let $f:F_1\rightarrow F_2$ be a  homomorphism. There exists a unique homomorphism $\phi:F_1/\textnormal{Ker~}f\rightarrow F_2$ such that the following diagram is commutative

	\begin{figure}[H]
			\centering
			\includegraphics[scale=.3]{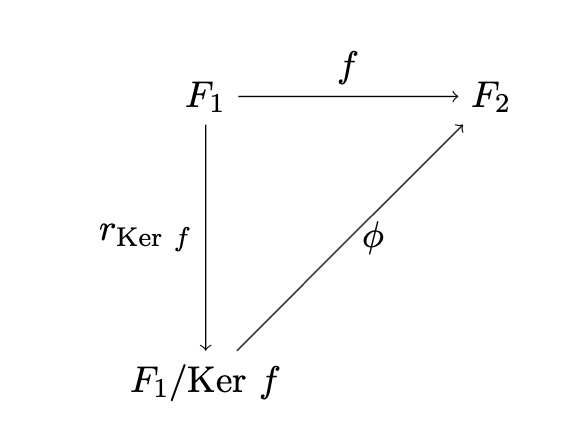}	
		\end{figure}
\end{theo}

\noindent that is, $\phi\circ r_{\textnormal{Ker~}f}=f$, where $r_{\textnormal{Ker~}f}(x):= \rho_{\textnormal{Ker~}f}(x)$. Moreover, $\phi$ is onto and it has a trivial kernel, namely, $\textnormal{Ker~}\phi=\{ \textnormal{Ker~}f  \}$. However, $\phi$ is an isomorphism if and only if $\rho_f=\rho_{\textnormal{Ker~}f } $.

\begin{theo}\label{thm58}Let $I$ and $J$ be $\Gamma$-order-ideals of a commutative monoid $T$ with $I\subseteq J$. Then
	$$\left(\bigslant{T}{I}\right)/\left(\bigslant{J}{I}\right)\cong T/J.$$
\end{theo}

\noindent \textbf{Proof:} Define $f:T/I\rightarrow T/J$ by $f(\rho_I(g))=\rho_J(g)$ for all $\rho_I(g)\in T/I$.\\
\indent Let $\rho_I(g_1),\rho_I(g_2)\in T/I$ and suppose 
$\rho_I(g_1)=\rho_I(g_2)$. Then $(g_1+I)\cap (g_2+I)\neq \varnothing.$ Thus, $g_1+w_1=g_2+w_2$ for some $w_1,w_2\in I\subseteq J$. Thus, 
$(g_1+J)\cap (g_2+J)\neq \varnothing.$
By Remark \ref{rem36}(ii), $\rho_J(g_1)=\rho_J(g_2)$. Thus, $f(\rho_I(g_1))=f(\rho_I(g_2))$. Hence, $f$ is well-defined.\\
\indent Let $\rho_I(g_1),\rho_I(g_2)\in T/I$. Then
\begin{eqnarray*}
	f(\rho_I(g_1)+\rho_I(g_2)) &=& f(\rho_I(g_1+g_2))\\
	&=& \rho_J(g_1+g_2)\\
	&=& \rho_J(g_1)+\rho_J(g_2)\\
	&=& f(\rho_I(g_1))+f(\rho_J(g_2)).
\end{eqnarray*}
Hence, $f$ is a homomorphism.\\
\indent Let $\rho_I(g)\in \textnormal{Ker}~f$. Then $f(\rho_I(g))=\rho_J(0)$, the identity in $T/J$. Thus, $\rho_J(g)=\rho_J(0)$. By Corollary \ref{cor36b}, $g\in J$. Hence, $\rho_I(g)\in J/I$. Thus, $\textnormal{Ker}~f\subseteq J/I$. Let $\rho_I(g)\in J/I$. Then $g\in J$. By Corollary \ref{cor36b}, $\rho_J(g)=\rho_J(0)$. Thus, $f(\rho_I(g))
=\rho_J(g)=\rho_J(0)$. Accordingly, $\rho_I(g)\in \textnormal{Ker~}f$. Hence, $J/I\subseteq \textnormal{Ker~}f$. So, $J/I= \textnormal{Ker~}f$.\\
\indent For $\rho_I(x), \rho_I(y)\in T/I$, recall that $\rho_I(x)~\rho_f~ \rho_I(y)$ if and only if $f(\rho_I(x))=f(\rho_I(y))$. We claim that  $\rho_f=\rho_{\textnormal{Ker~}f}$. Let $\rho_I(z)\in T/I$. We show that $\rho_f(\rho_I(z))=\rho_{\textnormal{Ker~}f}(\rho_I(z))$.\\
\indent Let $\rho_I(w)\in \rho_{\textnormal{Ker~}f}(\rho_I(z))$. Then 
$(\rho_I(z) +\textnormal{Ker~}f)\cap (\rho_I(w) +\textnormal{Ker~}f)\neq \varnothing.$
Thus, there exist $y_1,y_2\in \textnormal{Ker~}f$ such that
$\rho_I(z)+y_1=\rho_I(w)+y_2.$
Hence, 
$f(\rho_I(z))=f(\rho_I(z))+0=f(\rho_I(z))+f(y_1)=f(\rho_I(z)+y_1)$ and 
$f(\rho_I(w))=f(\rho_I(w))+0=f(\rho_I(w))+f(y_2)=f(\rho_I(w)+y_2).$
So, by well-definedness of $f$, we have $f(\rho_I(z))=  f(\rho_I(z)+y_1)= f(\rho_I(w)+y_2) =f(\rho_I(w))$. Accordingly, $\rho_I(w)\in \rho_f(\rho_I(z))$. Thus, $\rho_{\textnormal{Ker~}f}(\rho_I(z))\subseteq \rho_f(\rho_I(z))$.\\
\indent Now, let $\rho_I(w)\in \rho_f(\rho_I(z))$. Then $f(\rho_I(z))=f(\rho_I(w))$, that is, $\rho_J(z)=\rho_J(w)$. Thus, $(w+J)\cap (z+J)\neq \varnothing.$ This implies that there exist $h_1,h_2\in J$ such that $w+h_1=z+h_2$. Hence, $\rho_I(h_1), \rho_I(h_2)\in J/I=\textnormal{Ker~}f$. Consequently, 
$\rho_I(w)+\rho_I(h_1)=\rho_I(w+h_1)=\rho_I(z+h_2)=\rho_I(z)+\rho_I(h_2).$
This implies that 
$(\rho_I(w)+\textnormal{Ker~}f)\cap (\rho_I(z)+\textnormal{Ker~}f)\neq \varnothing.$
Hence, $\rho_I(w)\in \rho_{\textnormal{Ker~}f}(\rho_I(z))$. Accordingly, $\rho_f(\rho_I(z)  )\subseteq \rho_{\textnormal{Ker~}f}(\rho_I(z))$.\\
\indent Therefore, 
$\rho_f(\rho_I(z)  )= \rho_{\textnormal{Ker~}f}(\rho_I(z))$
for all $\rho_I(z)\in T/I$, that is, $\rho_f=\rho_{\textnormal{Ker~}f}$. By Theorem \ref{thm13},  these all imply that
\begin{center}
	~~~~~~~~~~~~~~~~~~~~~~~~$(\bigslant{T}{I})/(\bigslant{J}{I})=(\bigslant{T}{I})/ \textnormal{Ker~}f~\cong ~T/J.$ ~~~~~~~~~~~~~~~~~~~~~~~~\qed
\end{center}

Theorems \ref{thm1b}, Corollary \ref{cor2b}, and the Jordan-H\"{o}lder Theorem are monoid adaptations of the Baumslag's short proof in the group setting \cite{Baumslag}.

\begin{theo}\label{thm1b}Let $T$ be a refinement $\Gamma$-monoid and $Q$, $L$ and $N$ be $\Gamma$-order ideals of $T$ such that $L\subseteq Q$. Then
	$$\bigslant{Q}{(L+(Q\cap N))}\cong \bigslant{(Q+N)}{(L+N)}.$$
\end{theo}

\noindent \textbf{Proof:} Define $f:Q\rightarrow (Q+N)/(L+N)$ by $f(q)=\rho_{L+N}(q)$. Let $a,b \in Q$ and suppose that $a=b$. Then 
$\varnothing \neq (a+(L+N))\cap (a+(L+N))=(a+(L+N))\cap (b+(L+N))$.
Thus, $\rho_{L+N}(a)=\rho_{L+N}(b)$. Thus, $f(a)=f(b)$. Hence, $f$ is well-defined. Also, by the definition of the operation in the quotient monoid, we have $f(a+b)=\rho_{L+N}(a+b)=\rho_{L+N}(a)+\rho_{L+N}(b)=f(a)+f(b)$ which means $f$ is a homomorphism.\\
\indent Let $x\in \textnormal{Ker~}f$. Then $\rho_{L+N}(x)=f(x)=\rho_{L+N}(0)$. Thus, $x\in L+N$. Hence, for some $l\in L$ and $n\in N$, $l+n=x\in Q$. Since $Q$ is a $\Gamma$-order ideal, it follows that $n\in Q$. Hence 
$n\in Q\cap N$. Accordingly, $x=l+n\in L+(Q\cap N)$.

\indent Let $x\in L+(Q\cap N)$. Then $x=l+c$ for some $l\in L$ and $c\in Q\cap N\subseteq N$. Then $x=l+n\in L+N$ which implies $f(x)=\rho_{L+N}(x)=\rho_{L+N}(0)$, by Corollary \ref{cor36b}.  Thus, $L+(Q\cap N)\subseteq \textnormal{Ker~}f$. Consequently, $\textnormal{Ker~}f=L+(Q\cap N).$\\
\indent In order to use Theorem \ref{thm13} and conclude the isomorphism, we are left to show that $\rho_f=\rho_{\textnormal{Ker~}f}$, that is, $\rho_f(x)=\rho_{\textnormal{Ker~}f}(x)$ for all $x\in Q$. \\
\indent Let $y\in \rho_{\textnormal{Ker~}f}(x)$. Then $(x+\textnormal{Ker~}f)\cap (y+\textnormal{Ker~}f)\neq \varnothing$. Thus, there exist $a,b\in \textnormal{Ker~}f$ such that $x+a=y+b$. Now, 
$$f(x)= f(x)+\rho_{L+N}(0)=f(x)+f(a)=f(x+b)$$ and 
$$f(y)= f(y)+\rho_{L+N}(0)=f(y)+f(a)=f(y+b).$$
Accordingly, $f(x)=f(x+a)=f(y+b)=f(y)$, that is, $y\in \rho_f(x)$. Hence, $\rho_{\textnormal{Ker~}f}(x)\subseteq \rho_f(x)$.\\
\indent 
{ Let $y\in \rho_f(x)$. Then $y \in Q$ and $(x+(L+N))\cap (y+(L+N))\neq \varnothing $.\\ Hence, there exist $a,b \in (L+N)$ such that $x+a = y+b$. Since $T$ is a refinement monoid, we have $x=e_1+e_2$, $a=e_3+e_4$, $y=e_1+e_3$ and $b=e_2+e_4$ for some $e_1,e_2,e_3,e_4\in T$. Since $(L+N)$ and $Q$ are $\Gamma$-order ideals, we find $e_1+e_2, e_1+e_3 \in Q$ and $e_3+e_4, e_2+e_4 \in (L+N)$. 
	Hence $e_2, e_3 \in (L+N)\cap Q \subset \textnormal{Ker~}f$. Obviously
	$x+e_3 =e_1 +e_2+e_3 = y+e_2$ and $e_2, e_3 \in \textnormal{Ker~}f$. Thus, $(x+ \textnormal{Ker~}f)\cap (y+\textnormal{Ker~}f) \neq \varnothing$ which gives $y \in  \rho_{\textnormal{Ker~}f}(x).$ \\
}
\indent Consequently, $\rho_{\textnormal{Ker~}f}(x)= \rho_f(x)$ for every $x\in Q$, Thus, $\rho_{\textnormal{Ker~}f}= \rho_f$. By Theorem \ref{thm13}, we have 
$$~~~~~~~~~~~~~~~~~\bigslant{Q}{\textnormal{Ker~}f}=\bigslant{Q}{(L+(Q\cap N))}\cong \bigslant{(Q+N)}{(L+N)}.~~~~~~~~~~~~~~~~~~\qed$$

\begin{coro}\label{cor2b} Given $A, A', B$ and $B'$ being $\Gamma$-order ideals of a refinement $\Gamma$-monoid $T$ such that $A'\subseteq A$ and $B'\subseteq B$. Then
	$$\bigslant{(A\cap B)+B'}{(A'\cap B)+B'}\cong \bigslant{(A\cap B)+A'}{(A'\cap B)+A'}.$$
\end{coro}

\noindent \textbf{Proof:} By Theorem \ref{thm1b} with $N=B' $, $L=A'\cap B$, and $Q=A\cap B$, we obtain 
\begin{eqnarray*}
	\bigslant{A\cap B}{((A'\cap B)+(A\cap B\cap B'))}& =& \bigslant{Q}{(L+(Q\cap N))}\\
	&\cong& \bigslant{Q+N}{L+N}\\
	&=& \bigslant{(A\cap B)+B'}{(A'\cap B)+B'}.
\end{eqnarray*}
That is, 
$$\bigslant{A\cap B}{((A'\cap B)+(A\cap B'))}\cong \bigslant{(A\cap B)+B'}{(A'\cap B)+B'}.$$
Similarly, for $N=A'$, we also have 
$$\bigslant{A\cap B}{((A'\cap B)+(A\cap B'))}\cong \bigslant{(A\cap B)+A'}{(A'\cap B)+A'}.$$
Combining, we have
$$~~~~~\bigslant{(A\cap B)+B'}{(A'\cap B)+B'}\cong \bigslant{(A\cap B)+A'}{(A'\cap B)+A'}.~~~~~\qed$$

\begin{defn}\label{df39}\textnormal { $T$ is said to be a $\Gamma$-$\emph{Noetherian}$ monoid if for every chain
		$$A_1 \subseteq A_2 \subseteq A_3 \subseteq \cdots $$
		of $\Gamma$-order-ideals of $T$, there is an integer $n$ such that $A_i = A_n$ for all $i\geq n$.
	}
	
\end{defn}

\begin{defn}\label{df40}\textnormal { $T$ is said to be a  $\Gamma$-$\emph{Artinian}$  monoid if for every chain
		$$B_1 \supseteq B_2 \supseteq B_3 \supseteq \cdots $$
		of $\Gamma$-order-ideals of $T$, there is an integer $m$ such that $B_i = B_m$ for all $i\geq m$.
	}
	
\end{defn}

\begin{rema}\label{Rem4z} We have the following properties of $\Gamma$-order ideals inherited from $\Gamma$-monoids:
\begin{enumerate}
\item[\textnormal{(i)}]
A $\Gamma$-order ideal of a $\Gamma$-Noetherian $\Gamma$-monoid is $\Gamma$-Noetherian. 
\item[\textnormal{(ii)}] A $\Gamma$-order ideal of a $\Gamma$-Artinian $\Gamma$-monoid is $\Gamma$-Artinian.  
\end{enumerate}
\end{rema}


\begin{defn}\label{dfh1}Let $I$ be a $\Gamma$-order ideal of $T$. We say
	\begin{enumerate}
		\item[\textnormal{i)}] $I$ is a $cyclic$ $ideal$ if for any $x\in I$, there is an $\alpha\in 
		\Gamma$ such that $^\alpha x=x$;
		\item[\textnormal{ii)}] $I$ is a $comparable$ $ideal$ if for any $x\in I$, there is an $\alpha \in \Gamma$ such that $^\alpha x> x$;
		\item[\textnormal{iii)}] $I$ is a $non$-$comparable$ $ideal $ if for any $x\in I$, and any $\alpha \in \Gamma$, we have $^\alpha x\parallel x$.
	\end{enumerate}
	
\end{defn}

\begin{exa} Let $T=\N \oplus \N \oplus \N \oplus \N $ be a free abelian monoid with the action of $\Z$ on $T$ defined by $^1(a,b,c,d)=(d,a,b,c)$ and extended to $\Z$. Then $T$ is a cyclic monoid as for any $x\in T$ we have $^4x=x$. 
\end{exa}

\begin{defn}\label{df45}\textnormal{Let $T$ be a $\Gamma$-order-ideal. A $\Gamma$-series for $T$ is a sequence of $\Gamma$-order ideals $$~~~~~~~~~~~~~~~~~~~0=I_0 \subseteq I_1 \subseteq I_2 \subseteq \cdots \subseteq I_n=T.~~~~~~~~~~~~~~~~~~~(*) $$  The \emph{length} of a $\Gamma$-series is the number of its proper inclusions. A $refinement$ of $(*)$ is any $\Gamma$-series of the form 
		$$~~~~~~~~~~~~~~~~~~~0=I_0 \subseteq I_1 \subseteq \cdots \subseteq I_i\subseteq N\subseteq I_{i+1}\subseteq \cdots \subseteq  I_n=T,~~~~~~~~~~~~~~~~~~~ $$ and this refinement is said to be $proper$ if $I_i\subsetneq N\subsetneq I_{i+1}$.
		Furthermore, we say $(*)$ is a  $\Gamma$-$composition$ $series$ if for each  $i=0,1\cdots,n-1$, $I_i\subsetneq I_{i+1}$ and each of  quotients $I_{i+1}/I_{i}$ are simple $\Gamma$-monoids.\\ \indent  A $\Gamma$-composition series is of $cyclic$ (\textit{non-comparable}, \textit{comparable}) type if all of the simple quotients $I_{i+1}/I_{i}$ are cyclic (non-comparable, comparable). We further say a composition series is of $mixed$ $type$ of certain kinds if the simple quotients are those given kinds. We say two composition series are $equivalent$ if there is a one-to-one correspondence between the simple quotients of the composition series such that the corresponding quotients are $\Gamma$-isomorphic monoids.
	}
\end{defn}

\begin{exa} Consider the free abelian monoid
	\begin{multline*}
		T=\big\langle  a(i), b(i), c(i): a(i)=a(i+1)+b(i+1),\\ b(i)=b(i+1)+c(i+1), c(i)=c(i+1), i\in \Z\big\rangle.
	\end{multline*}
	It could be shown that $T$ is a $\Z$-monoid under the action $^nv(i):=v(i+n)$. $0\subseteq\langle c(0)\rangle \subseteq T$ is a $\Z$-series for $T$ which has a proper refinement $0\subseteq\langle c(0)\rangle  \subseteq \langle c(0), b(0)\rangle \subseteq T$ which is actually a $\Z$-composition series.
\end{exa}

\begin{theo}\textnormal{(Jordan-H\"{o}lder)} \label{thmh2}
	Two $\Gamma$-series of a refinement  $\Gamma$-monoid $T$ have equivalent refinement. Thus, any $\Gamma$-composition series are equivalent and a $\Gamma$-monoid having a composition series determines a unique list of simple $\Gamma$-monoids.
\end{theo}

\noindent \textbf{Proof:} Let 
$$ ~~~~~~~~~~~~~~~~ T= G_0 \supseteq G_1\supseteq G_2 \supseteq \cdots \supseteq G_n ~~~~~~~~~~~~~~~~~(*)$$
and $$ ~~~~~~~~~~~~~~~~ T= H_0 \supseteq H_1\supseteq H_2 \supseteq \cdots \supseteq H_m ~~~~~~~~~~~~~~~~~(**)$$
be two $\Gamma$-series for $T$. Let $G_{n+1}=\{0\}=H_{m+1}$. Now, for each $i=0,1,\cdots , n-1$, we consider the $\Gamma$-order ideals 
$$G_{i+1}+(G_i\cap H_j),$$
$j=0,1,\cdots, m+1$. Now, by Lemma \ref{lemc3},   $G_{i+1}+(G_i\cap H_j)$ is a $\Gamma$-order-ideal of $T$. Now, we consider the chain of $\Gamma$-order-ideals
\begin{eqnarray*}
	G_i&=&G_{i+1}+(G_i\cap H_0)\supseteq G_{i+1}+(G_i\cap H_1)\supseteq \cdots \supseteq G_{i+1}+(G_i\cap H_{m})\\
	&~& ~~\supseteq G_{i+1}+(G_i\cap H_{m+1}) = G_{i+1}. 
\end{eqnarray*}
Denote $G_{i+1}+(G_i\cap H_{j})$ by $G(i,j)$ and we obtain a refinement for $(*)$
\begin{eqnarray*}
	T&=& G(0,0)\supseteq G(0,1)\supseteq \cdots \supseteq  G(0,m)\subseteq G(0,1)\supseteq G(1,1)\\ &~& \supseteq \cdots \supseteq G(1,m)\supseteq G(n-1,m)\supseteq G(n,0)\supseteq G(n,1)\\ &~&~~~ \supseteq G(n,m)\supseteq G(n,m+1) \hspace{5cm}(*')
\end{eqnarray*}
Notice that $(*')$ has $(n+1)(m+1)$ (not necessarily distinct) terms.\\
\indent Similarly, for $(**)$, we obtain a refinement:
\begin{eqnarray*}
	T&=& H(0,0)\supseteq H(1,0)\supseteq \cdots \supseteq H(n,0)\subseteq H(0,1)\supseteq H(1,1)\\ &~& \supseteq \cdots \supseteq H(n,1) \supseteq H(0,2)\supseteq H(1,2)\supseteq\cdots \supseteq H(n-1,m)\\ &~&~~~ \supseteq H(n,m)\supseteq H(n+1,m)\hspace{5cm}(**')
\end{eqnarray*} which also has $(n+1)(m+1)$ terms. \\
\indent Now by Corollary $\ref{cor2b}$, we have

$$\dfrac{G(i,j)}{G(i,j+1)}= \dfrac{G_{i+1}(G_i\cap H_j)}{G_{i+1}(G_i\cap H_{j+1})}\cong \dfrac{H_{j+1}(G_i\cap H_j)}{H_{i+1}(G_i\cap H_{j})}= \dfrac{H_(i,j)}{H(i+1,j)}.$$
Thus, $(*')$ and $(**')$ are equivalent. \qed

\begin{rema}\label{Rem1z}
	If a composition series $\alpha$ exists for a $\Gamma$-monoid $T$, then the length of any $\Gamma$-series of $T$ is at most the length of $\alpha$.
\end{rema}

\begin{defn}\textnormal{
		A $\Gamma$-monoid $T$ is said to satisfy \emph{maximal condition} if every nonempty set of $\Gamma$-order ideals of $T$ has a maximal element under set-theoretic inclusion. }
\end{defn}

\begin{theo}\label{thm2z}
	For a $\Gamma$-monoid $T$, the following are equivalent.
	\begin{enumerate}
		\item [\textnormal{(i)}] $T$ is $\Gamma$-Noetherian
		\item [\textnormal{(ii)}] $T$ satisfies maximal condition
	\end{enumerate}
\end{theo}
\noindent \textbf{Proof:} (i) $\Rightarrow$ (ii) Suppose $T$ is $\Gamma$-Noetherian and let $\Sigma$ be a nonempty set of $\Gamma$-order ideal of $T$. Assume in the contrary that $\Sigma$ has no maximal element. Since $\Sigma\neq \varnothing$, there exists a $\Gamma$-order ideal $M_1\in \Sigma$. Since $\Sigma$ has no maximal element. there exists $M_2\in \Sigma$ such that $M_1\subsetneq M_2.$
Again, $M_2$ is not maximal in $\Sigma$. Hence there exists $M_3\in \Sigma$ such that $M_1\subsetneq M_2 \subsetneq M_3.$
Continuing in this manner, we have an ascending chain of $\Gamma$-order ideals in $\Sigma$
$M_1\subsetneq M_2 \subsetneq M_3\subsetneq \cdots .$ 
This contradicts our assumption that $T$ us $\Gamma$-Noetherian. Hence, $T$ must satisfy the maximal condition. \\
(ii) $\Rightarrow$ (i) Let 
$$~~~~~~~~M_1\subsetneq M_2 \subsetneq M_3\subsetneq \cdots ~~~~~~~~(*)$$
be any ascending chain of $\Gamma$-order ideals of $T$ and let 
$$\Sigma=\{M_i:i=1,2,3,\cdots  \}.$$
Then $\Sigma\neq \varnothing$. By our assumption, $\Sigma$ has a maximal element (say) $M_n$ for some $n\in \N$. In view of $(*)$, we have 
\begin{center}
	$M_n\subseteq M_k$
	for all $k>n$.
\end{center}
Since $M_n$ is a maximal element of $\Sigma$, we must have \begin{center}
	$M_n=M_k$
	for all $k>n$.
\end{center}
Hence, $T$ is $\Gamma$-Noetherian. \hfill \qed

\begin{theo}\label{thm3z}
	The following statements are equivalent.
	\begin{enumerate}
	\item[\textnormal{(i)}] $T$ has a $\Gamma$-composition series.
	\item[\textnormal{(ii)}] $T$ is $\Gamma$-Noetherian and $\Gamma$-Artinian. 
	\end{enumerate}
\end{theo}
\noindent \textbf{Proof:} Suppose $T$ has a $\Gamma$-composition series of length $n$. Then, by the Jordan-H\"{o}lder Theorem, it follows that every $\Gamma$-series has length at most $n$.\\
\indent Suppose $T$ is not $\Gamma$-Noetherian. Then there exists a chain of $\Gamma$-order ideals
$$N_1\subsetneq N_2\subsetneq N_3\subsetneq \cdots$$
for $T$ such that $N_i\neq N_{i+1}$ for all $i\in \N$. Hence, we have a $\Gamma$-series
$$\{0\}=N_0\subsetneq N_1\subsetneq N_2\subsetneq \cdots N_n\subsetneq N_{n+1}\subsetneq T$$
of length $n+1$. By Remark \ref{Rem1z}, we obtain a contradiction. Hence, $T$ is $\Gamma$-Noetherian. \\
\indent Suppose $T$ is not $\Gamma$-Artinian. Then there exists a chain of $\Gamma$-order ideals
$$N_1\supsetneq N_2\supsetneq N_3\supsetneq \cdots$$
for $T$ such that $N_i\neq N_{i+1}$ for all $i\in \N$. Hence, we have a $\Gamma$-series
$$T=N_0\supsetneq N_1\supsetneq N_2\supsetneq \cdots N_n\supsetneq N_{n+1}\supsetneq \{0\}$$
of length $n+1$. By Remark \ref{Rem1z}, we obtain a contradiction. Hence, $T$ is $\Gamma$-Artinian.\\
\indent Conversely, suppose $T$ is both $\Gamma$-Noetherian and $\Gamma$-Artinian. We show that $T$ has a $\Gamma$-composition series. If $T$ is simple, then we are done. Suppose $T$ is not simple and let $\Sigma_0$ be the set of all proper $\Gamma$-order ideals of $T$. Since $T$ is $\Gamma$-Noetherian, by Theorem \ref{thm2z}, $\Sigma_0$ has a maximal element, say $M_1$, that is, $M_1$ is a maximal $\Gamma$-order ideal of $T$. Now, if $M_1=\{0\}$, then we have a $\Gamma$-series $$~~~~~~T=M_0\supseteq M_1=\{0\}~~~~~~(*)$$
of length $1$. Since $M_1$ is maximal, $T/M_1$ is simple by Lemma \ref{lem57}. Thus, $(*)$ is a $\Gamma$-composition series for $T$. Suppose $M_1\neq \{0\}$. Then by Remark \ref{Rem4z}, $M_1$ is also Noetherian. Let $\Sigma_1$ be the set of all proper $\Gamma$-order ideals of $M_1$. Similarly, we obtain a maximal $\Gamma$-order ideal $M_2$ of $M_1$ and if $M_2=\{0\}$, we have a $\Gamma$-composition series
$$T=M_0\supsetneq M_1\supsetneq M_2=\{0\}$$
of length $2$. If $M_2\neq \{0\}$, we continue in the same manner ans obtain a strictly descending chain
$$T=M_0\supsetneq M_1\supsetneq M_2=\supsetneq \cdots \supsetneq M_i\supsetneq M_{i+1}\supsetneq \cdots$$
of $\Gamma$-order ideals of $T$ such that $M_i/M_{i+1}$ is simple for all $i=0,1,2,\cdots$. Since $T$ is $\Gamma$-Artinian, there exists $m\in \N$ such that $M_k=M_m$ for all $k>m$. This implies that $M_m$ has no proper $\Gamma$-order ideal. Thus, we obtain a $\Gamma$-composition series for $T$:
$$~~~~~T=M_0\supsetneq M_1\supsetneq M_2=\supsetneq \cdots \supsetneq M_m\supsetneq \{0\}.
$$\qed

\begin{lemm}\label{lem59}Let  $I$ be a $\Gamma$-order-ideal of $T$. Then $T$ has a composition series if and only if $T/I$ and $I$ have composition series.\end{lemm}

\noindent \textbf{Proof:}  Let 
$$0=I_n\subsetneq I_{n-1}\subsetneq \cdots \subsetneq I_1\subsetneq I_0=I$$
and 
$$0=\bigslant{T_m}{I}\subsetneq \bigslant{T_{m-1}}{I}\subsetneq \cdots \subsetneq \bigslant{T_0}{I}=\bigslant{T}{I}$$
be $\Gamma$-composition series for $I$ and $T/I$, respectively. Then
\begin{center}
	$\bigslant{I_{i}}{I_{i+1}}$ and $(\bigslant{T_{j}}{I})/(\bigslant{T_{j+1}}{I})\cong \bigslant{T_{j}}{T_{j+1}}$
\end{center}
are simple. Since $T_m/I=0=\{ \rho_I(0) \}$, $T_m=I$. Consider 
$$\hspace{1cm}0\subsetneq I_n\subsetneq \cdots \subsetneq I_0=I=T_m \subsetneq T_{m-1} \subsetneq \cdots \subsetneq T_1 \subsetneq T_0=T.\hspace{0.5cm}(*)$$
Then each of the factor is simple. Thus, $(*)$ is a $\Gamma$-composition series for $T$.

Conversely, suppose $T$ has a $\Gamma$-composition series
$$0=J_n\subsetneq J_{n-1}\subsetneq \cdots \subsetneq J_1\subsetneq J_0=T.$$
Then each quotient $J_i/J_{i+1}$ is simple. For each $i=1,2,...,n$, consider the $\Gamma$-order-ideal $I\cap J_i$. Then $(I\cap J_i)/(I\cap J_{i+1})$ is simple for each $i=1,2,...,n$. Hence, 
$$0=I\cap J_n\subsetneq I\cap J_{n-1}\subsetneq \cdots \subsetneq  I\cap J_1\subsetneq I\cap J_0=I\cap T=I$$
is a $\Gamma$-composition series for $I$.\qed

~\\
The following corollary directly follows from the proof of Lemma \ref{lem59}.

\begin{coro}\label{corg4}Let $I$ be a $\Gamma$-order ideal of $T$. If $I$ and $T/I$ have cyclic $[$resp., comparable, noncomparable$]$ composition series, then $T$ has cyclic $[$resp., comparable, noncomparable$]$ composition series.  \end{coro}

\begin{theo}\label{lemh3}Let $I_1, I_2,\cdots , I_k$ be distinct minimal $\Gamma$-order ideals of a refinement $\Gamma$-monoid $T$. Then
	$$ 0\subsetneq I_1\subsetneq I_1+I_2\subsetneq \cdots \subsetneq I_1+I_2+\cdots +I_k $$
	is a composition series for the $\Gamma$-monoid $I_1+I_2+\cdots +I_k $.
\end{theo}

\noindent \textbf{Proof:} Since $I_i$ are distinct $\Gamma$-order ideals, it is easy to show that the chain is proper. For $1<i\leq k$, the map
\begin{eqnarray*} \hspace{3cm} I_i& \longrightarrow & J_i \hspace{3cm} (*)\\
	 x & \mapsto & \rho_{J_i}(x)
\end{eqnarray*}
where  $$J_i=\dfrac{I_1+I_2+\cdots +I_i }{I_1+I_2+\cdots +I_{i-1} }$$ is clearly a surjective homomorphism of $\Gamma$-monoids. If $\rho_{J_i}(x)=\rho_{J_i}(y)$, then $x+a=y+b$, where $x,y\in I_i$ and $a,b\in I_1+I_2+\cdots +I_{i-1}$. Since $T$ is refinement, there are $z_1,z_2,z_3,z_4\in T$ such that $x=z_1+z_2$, $y=z_3+z_4$ and $a=z_1+z_3$, $b=z_2+z_4$. It follows that $z_2\in I_i\cap (I_1+I_2+\cdots +I_{i-1} )$ which leads to a contradiction. Thus, $(*)$ is an isomorphism of monoids, implying the quotients are simple. This proves the lemma. \hfill \qed

\begin{coro}\label{corg2} Let $T$ be a refinement $\Gamma$-monoid. 
Let $I_1, I_2,\cdots , I_k\subseteq T$ be distinct cyclic \textnormal{[}respectively, comparable, noncomparable\textnormal{]} minimal $\Gamma$-order ideals. Then
	$$0\subsetneq I_1\subsetneq I_1+I_2\subsetneq \cdots \subsetneq I_1+I_2+\cdots +I_k \qquad (*)$$
	is a cyclic \textnormal{[}comparable, noncomparable\textnormal{]} $\Gamma$-composition series for the monoid $I_1+I_2+\cdots +I_k $.
\end{coro}
\noindent \textbf{Proof:} By Theorem \ref{lemh3}, we are left to show that the quotients 
$$\bigslant{I_1+I_2+\cdots +I_t }{I_1+I_2+\cdots +I_{t-1} }$$
is cyclic for each $t=1,2,\cdots , k$. Note that by Lemmas \ref{lemc3} and \ref{lem46}, $I_1+I_2+\cdots +I_{t-1}\cap I_t $ is a $\Gamma$-order ideal of $I_t$.  Now, since the $I_t$ is minimal, it follows that    $I_1+I_2+\cdots +I_{t-1}\cap I_t=\{0\} $. Hence, by Lemma \ref{lemg1}, we have
$$\bigslant{I_1+I_2+\cdots +I_t }{I_1+I_2+\cdots +I_{t-1} }\cong I_t,$$
which is cyclic. Hence, $(*)$ is a cyclic composition series for $I_1+\cdots+I_k$. The same follows if the $I_i's$ are comparable [resp., noncomparable] ideals. \hfill \qed

\section*{Acknowlegdement} The Ph.D. scholarship of A.~Sebandal within the Accelerated Science and Technology Human Resource Development Program (ASTHRDP) under the Department of Science and Technology - Science Education Institute (DOST-SEI) is gratefully acknowledged. We want to emphasize that without this funding, this paper would not be able to be finished.

\section*{Conflict of Interest} None of the authors has a conflict of interest in the conceptualization or publication of this work.

\bibliographystyle{plain}

\begin{thebibliography}{20}
	\baselineskip=1\baselineskip
	
	
	\bibitem{Baumslag} Baumslag B. \textit{A simple way of proving the Jordan--H\"older--Schreier theorem.} The American Mathematical Monthly. 113(10), 933-935. (2006)
	
	
	\bibitem{Cassidy} P. J. Cassidy, and M. F. Singer, \textit{A Jordan-–H\"older theorem for differential algebraic groups.} Journal of Algebra, 328(1), 190-217. (2011).
	
	
	\bibitem{Erdmann} K. Erdmann and T. Holm. \textit{Algebras and representation theory.}  Springer Berlin, 2018.
	
	\bibitem{Fritsch} R. Fritsch. \textit{Variations on the theorems of Jordan-H\"older and Schreier.} Journal of pure and applied algebra.(2),  209-29. (1972).
	
	\bibitem{monoid.q} Y. Give'on, \emph{Normal monoids and factor monoids of commutative monoids}, The University of Michigan, 1963.
	
	
	\bibitem{A1}
	R. Hazrat, H. Li, \emph{The talented monoid of a Leavitt path algebra}, Journal of Algebra 547, 430-455. (2020)
	
	\bibitem{Lederman} P.J. Hilton and W. Ledermann \textit{On the Jordan‐H\"older Theorem in Homological Monoids.} Proceedings of the London Mathematical Society. 3(1), 321-34. (1960).
	
	
	\bibitem{Hol89} O. H\"older, \textit{Zur\"uckf\"uhrung einer beliebigen algebraischen Gleichung auf eine Kette von Gleichungen}. Mathematische Annalen. 34(1), 26-56. (1889).
	
	\bibitem{Jor69} C.M. Jordan,  \textit{Commentaire sur Galois.} Mathematische Annalen 1.2, 141-160. (1869).
	
	\bibitem{Jor70}
	C.M. Jordan \textit{Traité des substitutions et des équations algébriques}. Gauthier-Villars; 1870.
	
	
	\bibitem{Kashiwara} M. Kashiwara and M.K. Shapira, Categories and sheaves. Grundlehren der Math. Wiss. 2006.
	
	
	
	\bibitem{Natale} S.~Natale, \textit{Jordan-H\"older theorem for finite dimensional Hopf algebras}. Proceedings of the American Mathematical Society. 143(12), 5195-211. (2015).
	
	
	\bibitem{Natale2} S.~Natale  \textit{A Jordan–-H\"older theorem for weakly group-theoretical fusion categories.} Mathematische Zeitschrift. 1,283(1-2):367-79. (2016).
	
	\bibitem{Orlik} S. Orlik, and M. Strauch. \textit{On Jordan-H\"older series of some locally analytic representations.} Journal of the American Mathematical Society 28, no. 1, 99-157. (2015)
	
	\bibitem{Book1} J. Perry, \emph{Notes on algebra}, University of Southern Mississippi, 2009.
	
	
	
	\bibitem{Sch28} O. Schreier \textit{\"uber den Jordan-H\"olderschen Satz}. Abhandlungen aus dem Mathematischen Seminar der Universit\"at Hamburg, 6(1),300-302. Springer-Verlag. (1928).
	
	
	
\end{thebibliography}

\end{document}